\documentclass[12pt,english]{article}
\usepackage[cp1251]{inputenc}
\usepackage[english]{babel}
\usepackage[T1]{fontenc}
\usepackage{amsthm, amsfonts,amsmath,amssymb,graphicx, mathrsfs,url}

\providecommand{\keywords}[1]{\textbf{Keywords:} #1}
\providecommand{\subjclass}[1]{\textbf{MSC 2010:} #1}

\def\R{{\mathbb R}}

\newcommand{\true}{\operatorname{true}\nolimits}
\newcommand{\false}{\operatorname{false}\nolimits}
\def\Lie{\mathop{\rm Lie\,}\nolimits}
\def\xSE{\mathop{\rm SE\,}\nolimits}
\newcommand{\vect}[1]{\left( \begin{array}{c} #1 \end{array} \right)}
\newcommand{\sn}{\operatorname{sn}\nolimits}

\newcommand{\dn}{\operatorname{dn}\nolimits}
\newcommand{\sign}{\operatorname{sign}\nolimits}
\def\tcut{t_{\operatorname{cut}}}

\begin{document}

\title{Controlling of a mobile robot with a trailer \\ and its nilpotent approximation\footnote{This is a preprint of the Work accepted for publication in <<Regular and Chaotic Dynamics>>}}
\author{Andrey A.~Ardentov \\ \small{aaa@pereslavl.ru}  \\ \small{Program Systems Institute of RAS} \\ \footnotesize{Pereslavl-Zalessky, Yaroslavl Region, 152020, Russia}}
\maketitle
\begin{abstract} 
The work studies a number of approaches to solving motion planning problem for a mobile robot with a trailer. Different control models of car-like robots are considered from the differential-geometric point of view. The same models can be also used for controlling a mobile robot with a trailer. However, in cases where the position of the trailer is of importance, i.e., when it is moving backward, a more complex approach should be applied. At the end of the article, such an approach, based on recent works in sub-Riemannian geometry, is described. It is applied to the problem of reparking a trailer and implemented in the algorithm for parking a mobile robot with a trailer.
\end{abstract}



\section*{Introduction}
Control problems for wheeled robots are usually described by nonholonomic systems. A nonholonomic system arises when the dimension of configuration space is greater than the dimension of control. Here and below, the configuration space represents possible positions of a wheeled robot, i.e., for a car-like  robot it can be expressed in the following way: 
$$M_c \subseteq \R^2 \times S^1,$$ 
where $\R^2$ corresponds to the reference point of the robot on a plane and $S^1$ corresponds to the orientation of the robot. For nonholonomic systems,  motion in some directions is infinitesimally prohibited, but it is locally and globally possible (through complex maneuvers of the system). 

There is a well-known problem about parking of a car which can not move in a direction perpendicular to the direction of motion of the wheels. In fact, every driver has faced this problem. Let us consider the simplest case where a car should be parked on an empty parking lot (a plane) with no other cars or static obstacles. It is the classical hypothesis of "rolling without slipping" that provides the kinematic model of the car. The control problem for this model is stated as follows: given a system of differential equations, fixed initial and final positions of the mobile robot and restrictions on one or two dimensional control, one should find a control law and the corresponding trajectory satisfying these conditions. Note that in order to control a real robot, the kinematic model should be derived from a dynamic one.

\section{Kinematic control models of a car-like mobile robot}\label{car}
From a driver's point of view, cars have two controls: the accelerator and the steering wheel. The position of a rear-wheel drive car can be defined by a vector $q = (x,y,\theta) \in M_c$, where $(x,y) \in \R^2$ is the midpoint of the rear wheels and $\theta \in S^1$ is the angle of the car orientation which coincides with the direction of the rear wheels. The direction of the front wheels is not fixed and corresponds to the steering control. It is possible that each pair of wheels can be reduced to one wheel only. Moreover, if we are not concerned with the direction of the front wheels, the control system for a car-like robot is equivalent to the system for a wheel or a skate~\cite{hadham}:
\begin{align}
\dot{x} &= u_1 \cos \theta, \nonumber \\
\dot{y} &= u_1 \sin \theta, \label{sys} \\
\dot{\theta} &= u_2,	\nonumber
\end{align}
where $u_1, u_2$ are respectively linear and angular velocities as controls. Note that the mechanical constraint for the steering angle of the front wheels of a car should be reduced to the constraint on $u_1, u_2$. Such a system looks like the kinematic model of a unicycle or a two-driving wheel mobile robot. However, the dynamical models of these systems are not identical. The main difference between kinematic models lies in the admissible control domains which should be obtained from dynamical models. Further in the text, a car is identified with a mobile robot.

The problem is to find controls $\big(u_1 (t), u_2 (t)\big)$ and the corresponding trajectory $q_c(t) = \big(x(t),y(t), \theta(t)\big)$ satisfying system~(\ref{sys}) with boundary conditions:
\begin{align}
q_c(0)=(x_0, y_0, \theta_0), \qquad q_c(t_1) = (x_1, y_1, \theta_1), \label{positions}
\end{align}
where $(x_0,y_0,\theta_0), (x_1, y_1, \theta_1) \in M_c$.

Since such a formulation of the problem has a continuum of solutions in the general case, a cost functional should be added and minimized in order to choose the optimal solution among the family. 

The fact that we can choose the coordinate system $(x,y)$ allows us to always associate the initial position with the origin: 
\begin{align}
(x_0, y_0, \theta_0)=(0,0,0). \label{origin}
\end{align}

Further, the paper considers a number of known approaches to solve the described above problem, which admits different restrictions on the control. 

\subsection{Dubins car} \label{secdub}
A.A.~Markov is one of the first mathematicians to work on such a problem~\cite{markov}. In 1887, he studied four problems with application to railroad practice. The first problem can be formulated geometrically in the following way. Given two points $A,B \in \R^2$ in a plane, a positive number $R$ and a direction $AC \in S^1$, the problem is to find the shortest smooth path $\gamma$ joining $A$ and $B$ s.t. $AC$ is the tangent direction of $\gamma$ at the point $A$ and the curvature of $\gamma$ is bounded everywhere by $\frac{1}{R}$. 

\begin{figure}[ht]
\centering
\includegraphics[width=0.9\linewidth]{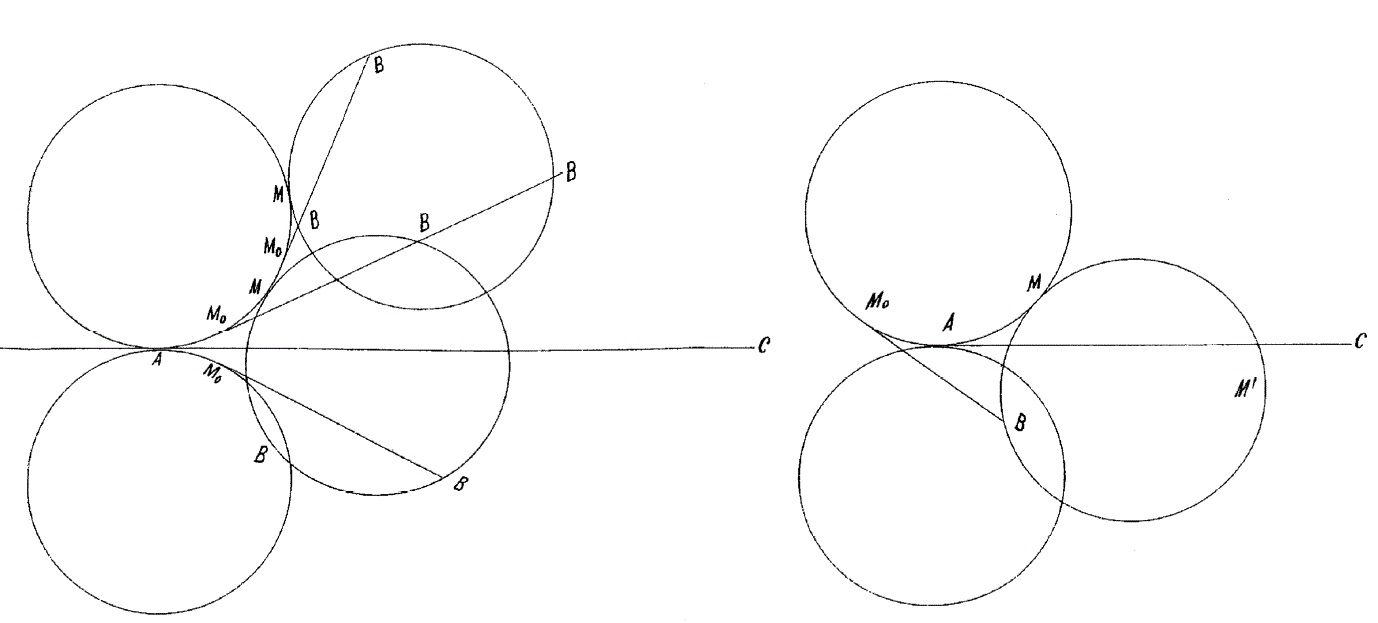}
\caption{Sketches from Markov's work~\cite{markov}}
\label{fig:skmark}
\end{figure}

The solution consists of an arc of a circle with radius $R$ and a straight segment or of two arcs depending on the position of the point $B$ (see Fig.~\ref{fig:skmark}). 
Such a solution can be applied to parking problem~(\ref{sys}),(\ref{positions}), provided that the car is moving with a  constant speed equal to one and is unable to move on a circle with radius smaller than $R$, i.e.,
\begin{align}
u_1 = 1, \qquad |u_2| \leq \frac{1}{R}, \label{dubins}
\end{align} 
and the direction of the car is of no importance at the desired end state, i.e., $\theta_1$ is not fixed. Since the car is moving with the unit velocity, the length minimization is equivalent to the time minimization: 
\begin{align}
t_1 \to \min. \label{timemin}
\end{align} 

Other problems from~\cite{markov} involve the condition about the direction at the point $B$, but also have some additional constraints on the curvature and its derivative. 

\begin{figure}[ht]
\centering
\includegraphics[width=0.7\linewidth]{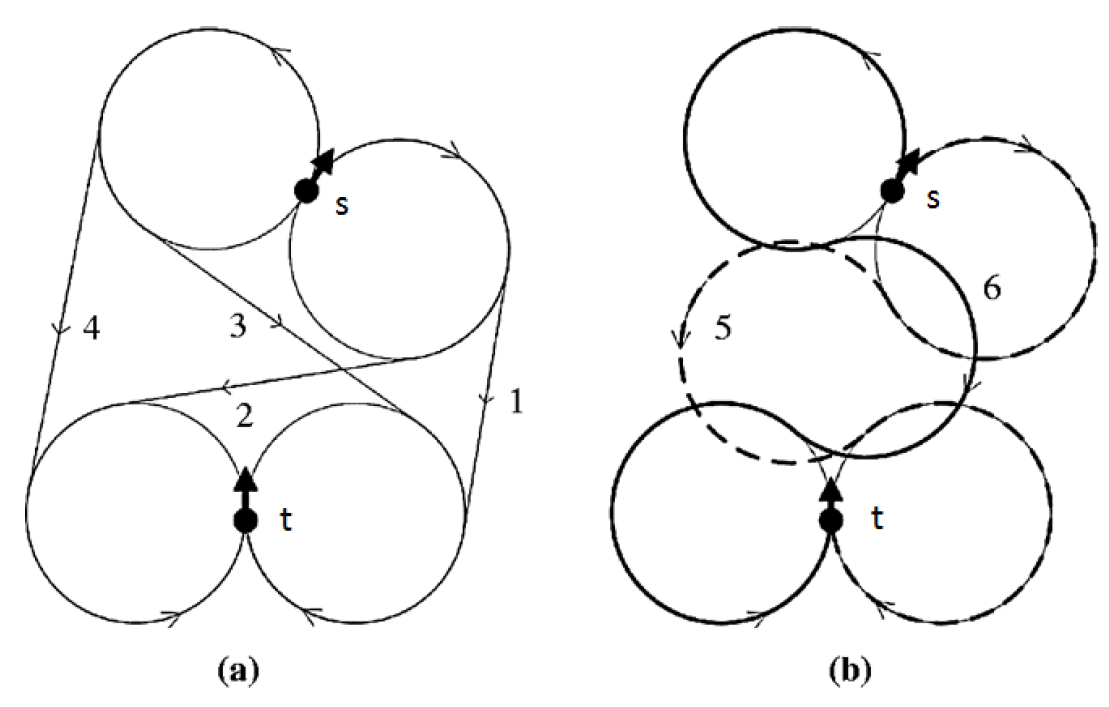}
\caption{6 variants of shortest paths for Dubins car}
\label{fig:dub}
\end{figure} 

In 1957,~L.~Dubins~\cite{dubins} solved problem~(\ref{sys})--(\ref{timemin}) with fixed $\theta_1$ using geometric methods. The solution is reduced to selection among 6 variants of smooth connection of arcs and straight segments (see Fig.~\ref{fig:dub}). This model corresponds to the time-optimal control problem for a car moving with a constant velocity and with bounds on angular velocity. Such a car is called the Dubins car. 

\subsection{Reeds-Shepp car}
In 1990, J.A.~Reeds and L.A.~Shepp~\cite{reeds} studied the so-called model of Reeds--Shepp car which can move forwards and backwards. This case gives the following constraints on controls:
\begin{align}
|u_1| = 1, \qquad |u_2| \leq \frac{1}{R}. \label{constraint2}
\end{align}

\begin{figure}[ht]
\centering
\includegraphics[width=0.3\linewidth]{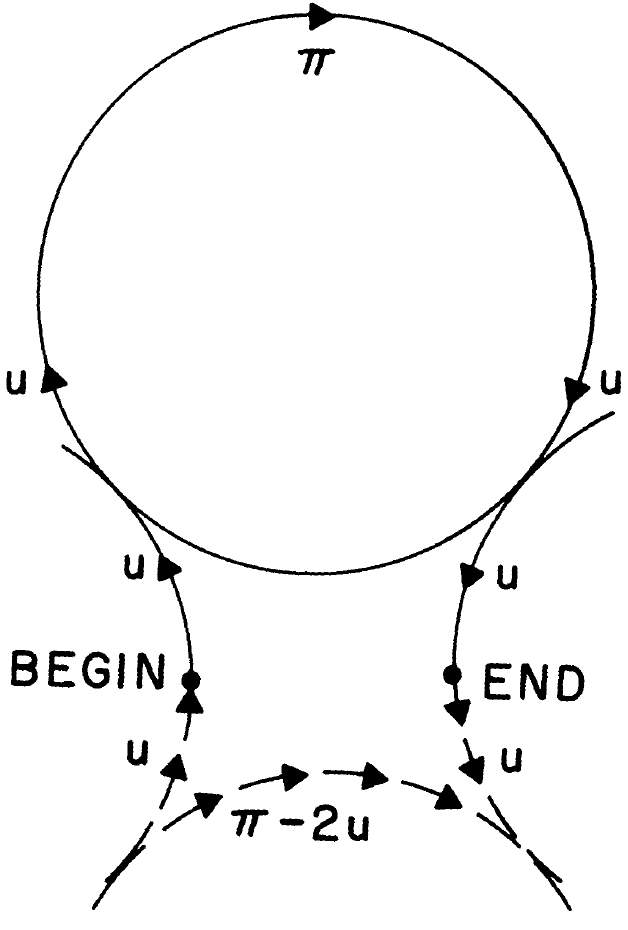}
\caption{Comparison of paths in Reeds-Shepp model~\cite{reeds}}
\label{fig:rsh}
\end{figure}

Such an admission allows the construction of paths shorter than in case~(\ref{dubins}). However, the number of variants of shortest paths increases significantly, the worst case providing 68 variants. Optimal trajectories could contain cusp points, i.e., points where the movement vector of the car changes to the opposite (see example in Fig.~\ref{fig:rsh}). As in Dubins' work, this characterization is done without obstacles. 

Soon, simpler solutions for both models were obtained in~\cite{sussmann} and in~\cite{boissonnat} via Pontryagin's maximum principle~\cite{PBGM}. In addition,  H.J.~Sussmann and G.~Tang reduced the sufficient family to 46 canonical paths. 

The Reeds--Shepp car is small time controllable from everywhere, i.e., the reachable set of admissible configurations for any small time $t>0$ contains a small neighborhood of an initial point (see geometric proof in~\cite{laum87}). The Dubins car is locally controllable, i.e., it is possible to find a path for the initial and final configurations that are arbitrarily close, but not small time controllable from everywhere, so the length of such a path does not converge to zero. Optimal control for both cars is a piecewise constant function. 

Both models serve as classical examples of time-optimal problems. However, another criterion for choosing an optimal solution which minimizes the amount of car maneuvers can be considered. A classic problem from mechanics with such a variational principle is described below.

\subsection{Planar elasticae} \label{secela}
\begin{figure}[ht]
\centering
\includegraphics[width=0.45\linewidth]{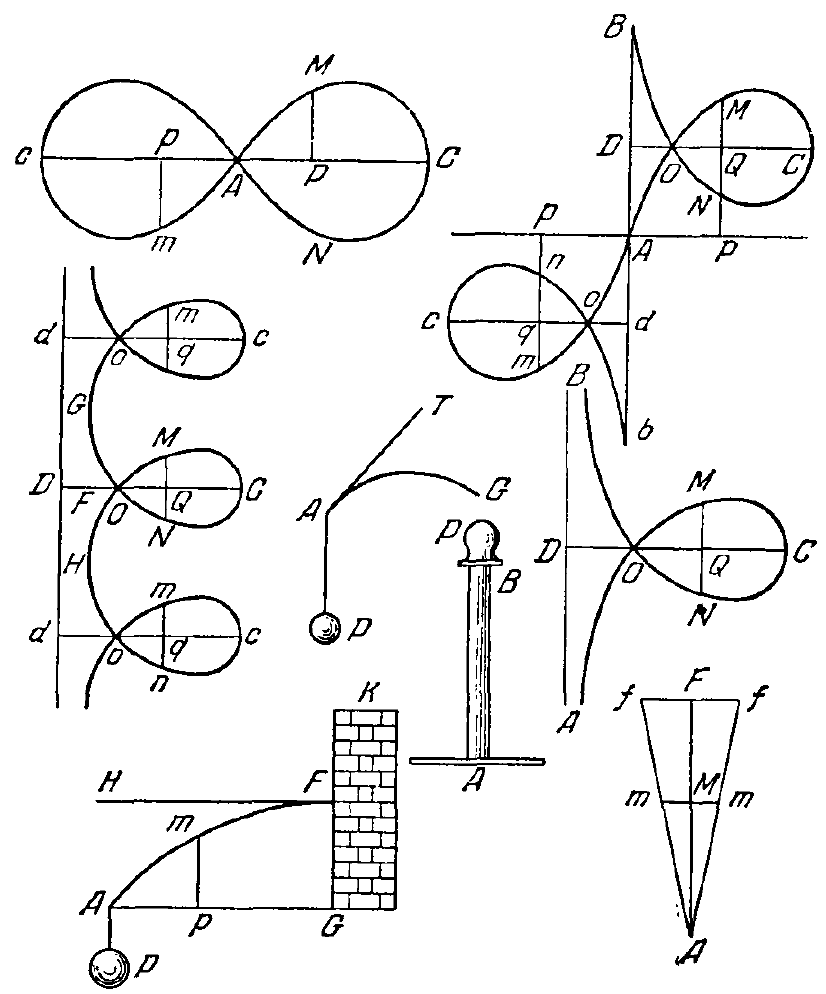} \includegraphics[width=0.47\linewidth]{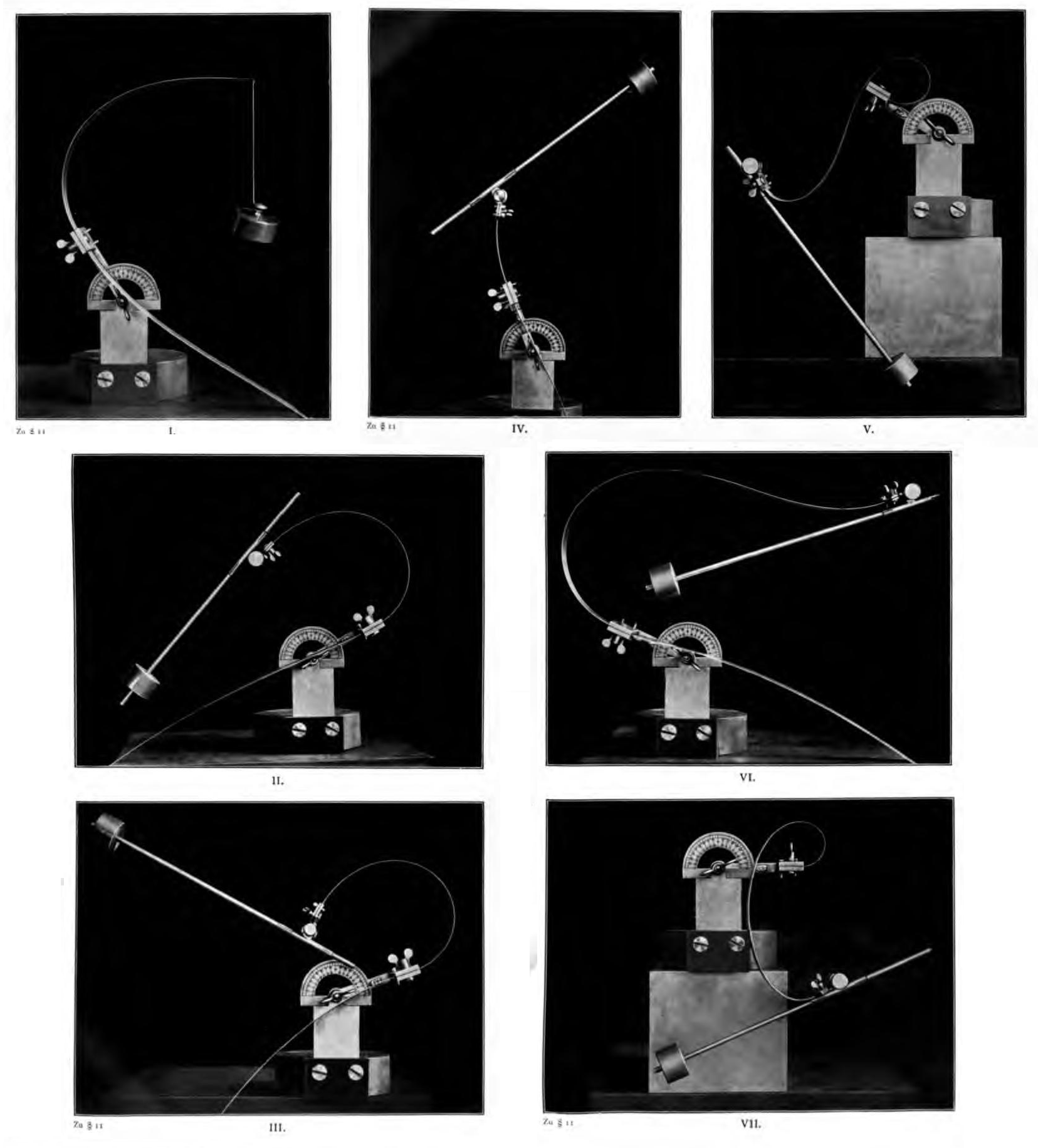}
\caption{Euler's scetches~\cite{euler} and Born's experiments~\cite{maxborn}}
\label{fig:euborn}
\end{figure}
In 1744, L.~Euler considered a problem on stationary configurations of an elastic rod~\cite{euler}. Given an elastic rod in a plane with fixed endpoints and tangents at the endpoints; should find the possible profiles of the rod with given boundary conditions. Euler derived differential equations for stationary configurations of a rod and described their possible qualitative types. These
configurations are called Euler's elasticae. In 1880, L.~Saalschutz obtained explicit parametrization of the curves~\cite{saal}. These curves have rich history including Max Born's work~\cite{maxborn} on stability of Euler's elasticae (see Fig.~\ref{fig:euborn}), for more details see~\cite{eus1}. They appear as solutions to several problems of sub-Riemannian geometry~\cite{notes}: nilpotent sub-Riemannian problem on the Engel group~\cite{engel}, nilpotent sub-Riemannian problem with the growth vector $(2,3,5)$~\cite{brockett}, the problem of rolling a sphere over a plane~\cite{jurd}.

In other words, we have system~(\ref{sys}) with $u_1 = 1, u_2 \in \R$, boundary conditions~(\ref{positions}) and a cost functional of energy to be minimized
\begin{align}
\int_0^{t_1} u_2^2 (t) \ d t \to \min. \label{inte}
\end{align}

One of the solutions to the parking problem of a car can be expressed via the shape of an elastic rod~\cite{walsh},  provided the car is moving forwards with a constant linear speed and with no boundary on steering. Time minimization in this case has no solutions (curves tend to straight segments but can not reach them for general boundary conditions). Note that time $t_1$ is fixed in~(\ref{inte}) and is equal to the length of the corresponding rod. Jacobi elliptic functions~\cite{jacobi} are used in the expression of optimal control~\cite{eus1} for a car moving along an elastica. 

The global structure of all solutions of the problem is described in a series of works~\cite{eus1, eus2, expela}, the software and algorithm for numerical computation of optimal elasticae are given in~\cite{solela}. A tool for plotting and evaluating generic Euler's elastica can be found in~\cite{demo}.

All models considered above have one-dimensional steering control $u_2$. The Reeds--Shepp car has an additional discrete control for direction of movement. A natural extension to it is the Reeds-Sheep car with varying speed and two dimensional control. One of such models arises in sub-Riemannian geometry~\cite{mont}.

\subsection{Differential drive robot} \label{uni}
Consider system~(\ref{sys}) with no bounds on control $(u_1, u_2) \in \R^2$. The problem is to find a curve connecting two points~(\ref{positions}) with a minimum sub-Riemannian length
\begin{align}
\int_0^{t_1} \sqrt{u_1^2(t) + u_2^2(t)} \ d t \to \min. \label{intsub}
\end{align}

Optimal solutions for this problem were obtained by Yu.~Sachkov in~\cite{se21,se22,se23} via geometric control theory~\cite{notes}.  A more general length functional reflecting the different weights of controls can be minimized:
\begin{align}
\int_0^{t_1} \sqrt{u_1^2(t) + \alpha^2 u_2^2(t)} \ d t \to \min, \qquad \alpha > 0. \label{intsub2}
\end{align}
The extension is easily implemented by changing the configuration variables $(x,y, \theta)$. The parameter $\alpha$ corresponds to the choice of scale in plane $(x,y)$.

This solution can be applied to a mobile robot with two parallel driving wheels, the acceleration of each being controlled by an independent motor. Such a robot is called a differential drive robot. The distance between the wheels may define the value of parameter $\alpha$.

\section{Geometric formalization of the control problem for a trailer-like robot}
\begin{figure}[ht]
\centering
\includegraphics[width=0.45\linewidth]{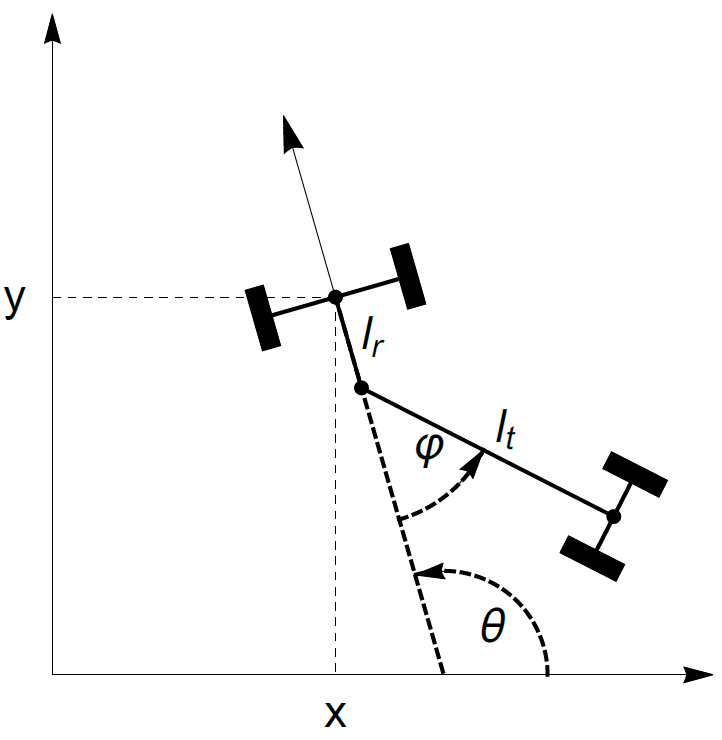}
\caption{Parameters of the system for a mobile robot with a trailer}
\label{fig:trailer}
\end{figure}

Consider a general model of a mobile robot with a  trailer (see Fig.~\ref{fig:trailer}). If a robot is able to move forward only then such a system can usually be reduced to system~(\ref{sys}), since the position of the trailer is of no importance. Such systems can use Dubins paths (see. Subsec.~\ref{secdub}) or planar elastica (see. Subsec.~\ref{secela}) for the path planning problem.

That being said, the complexity of the parking problem increases significantly for a car with a trailer moving backwards. A detailed survey on the various control strategies for the backward motion of a mobile robot with trailers can be found in~\cite{surv}. Formalization of control laws for such a problem requires the study of a nonlinear four-dimensional differential system with two-dimensional control which defines the linear and angular velocity. Such a differential system has the following form:
\begin{align}
\dot{q} &= u_1 X_1(q)  + u_2 X_2(q),  \label{sysgen} \\
q &= (x,y,\theta,\varphi)^T \in M=\R^2\times(S^1)^2, \quad (u_1, u_2) \in \R^2, \label{space}\\
X_1(q) &= \vect{\cos \theta  \\ \sin \theta \\  0 \\ - \frac{\sin \varphi}{l_t}}, \qquad X_2(q) = \vect{0 \\ 0 \\ 1 \\ - \frac {l_r \cos \varphi}{l_t} - 1} \label{vf},
\end{align}
where $x,y,\theta$ are coordinates of the position of the car defined in~Sec.~\ref{car}, $\varphi \in S^1$ is the angle of the trailer with respect to the car; constant parameters $l_r, l_t$ define distances from the robot and the trailer to the connection point as shown in Fig.~\ref{fig:trailer}. 

The problem is to move a robot with a trailer from one configuration to another, i.e., to find a path $q(t)$, s.t.
\begin{align}
q (0)= q_0=(x_0, y_0, \theta_0, \varphi_0), \qquad q (t_1) = q_1=(x_1, y_1, \theta_1, \varphi_1), \label{positions2}
\end{align}
where $q_0, q_1 \in M.$

Such a problem has symmetries which translate solutions of the problem to other solutions. The group of motions of the plane gives such symmetries and allows to move the initial position of a car to the origin~(\ref{origin}). Another known symmetry is the similarity transformation $\delta_{\mu}$, which changes configuration coordinates, constants and controls in the following way: 
\begin{align}
\delta_{\mu}: (x,y,\theta,\varphi, l_t, l_r, u_1, u_2) \mapsto (\mu x, \mu y, \theta, \varphi, \mu l_t, \mu l_r, \mu u_1, u_2), \label{sym}
\end{align}
this symmetry is used in the algorithm for trailer reparking, see Subsec.~\ref{reparking}.

Note that usually there is mechanical constraint $|\varphi| \leq \varphi_{\max}$. Such a constraint can be treated as an obstacle in $\R^2 \times (S^1)^2$.

\section{Obstacle avoidance}
The presence of obstacles can be described by inequalities and equalities for $x,y$. In other words we have a function $\R^2 \to \{\true, \false\}$ which tells us if a point $(x,y)$ belongs to an obstacle.

The resulting constraint equations for configuration space have different forms for different shapes of a robot. The simplest shape of a car-like robot is a circle because such constraint equations do not involve the angle parameter~$\theta$. The algorithm for computing the configuration space obstacles, when the objects are polygons, is given in~\cite{lozano}.

\begin{figure}[ht]
\begin{center}
\includegraphics[width=0.4\linewidth]{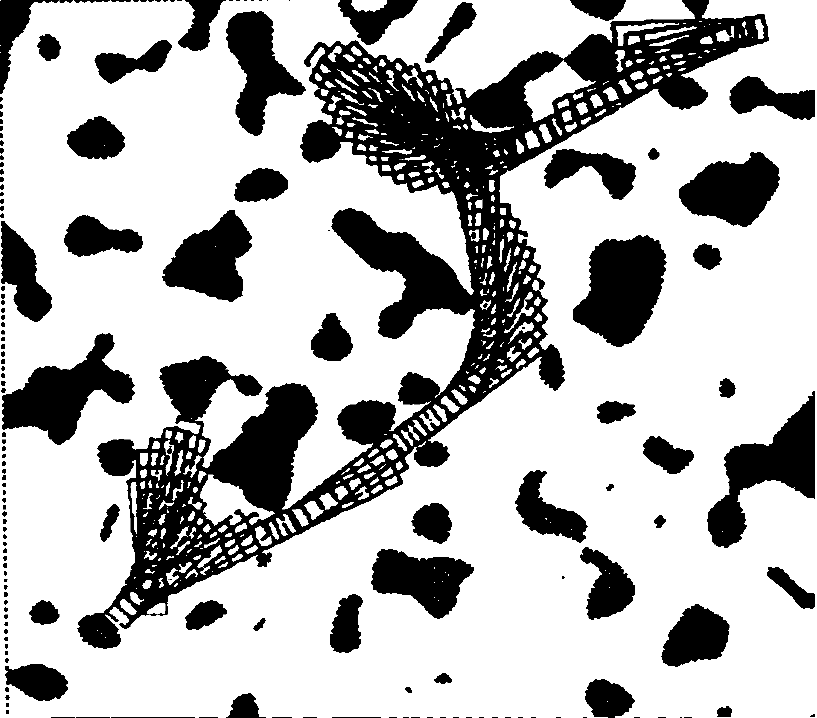} \qquad \includegraphics[width=0.4\linewidth]{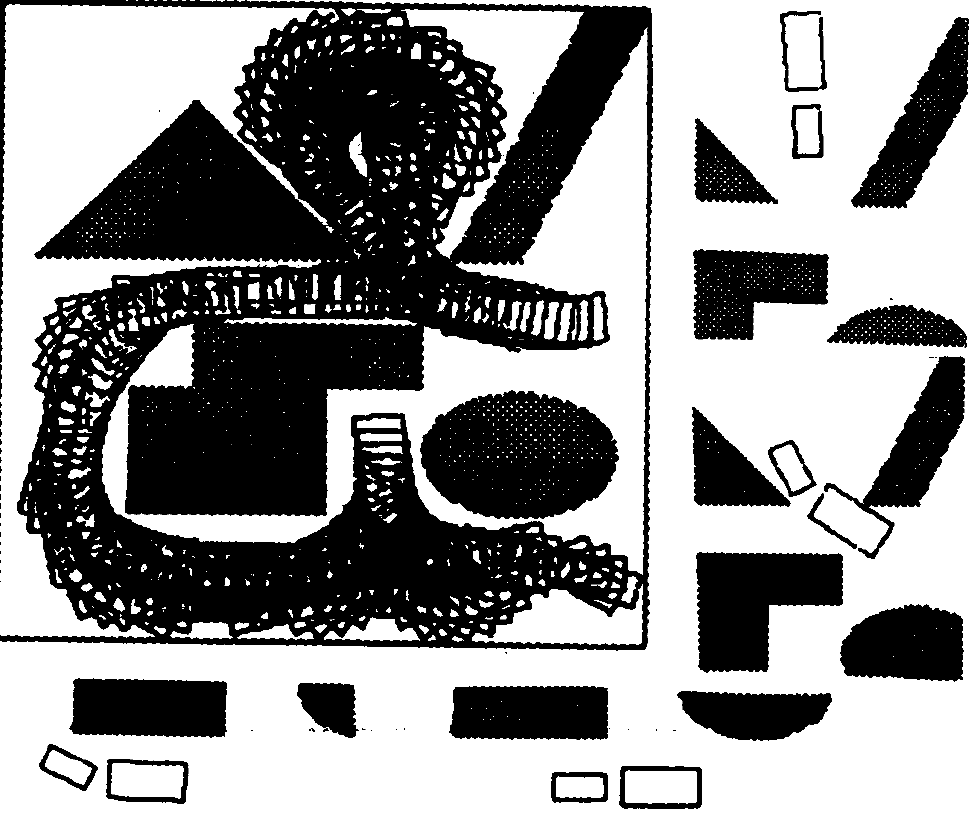}
\caption{Planner for car-like and trailer-like robots~\cite{barlat}}
\label{fig:parking}
\end{center}
\end{figure}

In 1989, J. Barraquand and J.-C. Latombe~\cite{barlat}, using standard results of differential geometry, nonlinear control theory and dynamic programming search, developed path planning algorithms for two robotic systems: car-like robots and trailer-like robots. A car-like robot is kinematically similar to a rectangular-shaped car. Left Fig.~\ref{fig:parking} shows an example of maneuvering in an unstructured workspace represented as a $512^2$ bitmap with a maximal steering angle of the front wheels equal to 45 degrees (the running time was about 2 minutes). A trailer-car is kinematically similar to a vehicle towing a trailer. Right Fig.~\ref{fig:parking} shows an example with a simulated two-body trailer-like robot where the trailer has to maneuver in a cluttered workspace with a maximal steering angle equal to 45 degrees (the running time was about 5 minutes). Their planners confirm controllability of such systems among obstacles in practice, but the solutions are not optimal.

\begin{figure}[ht]
\centering
\includegraphics[width=0.45\linewidth]{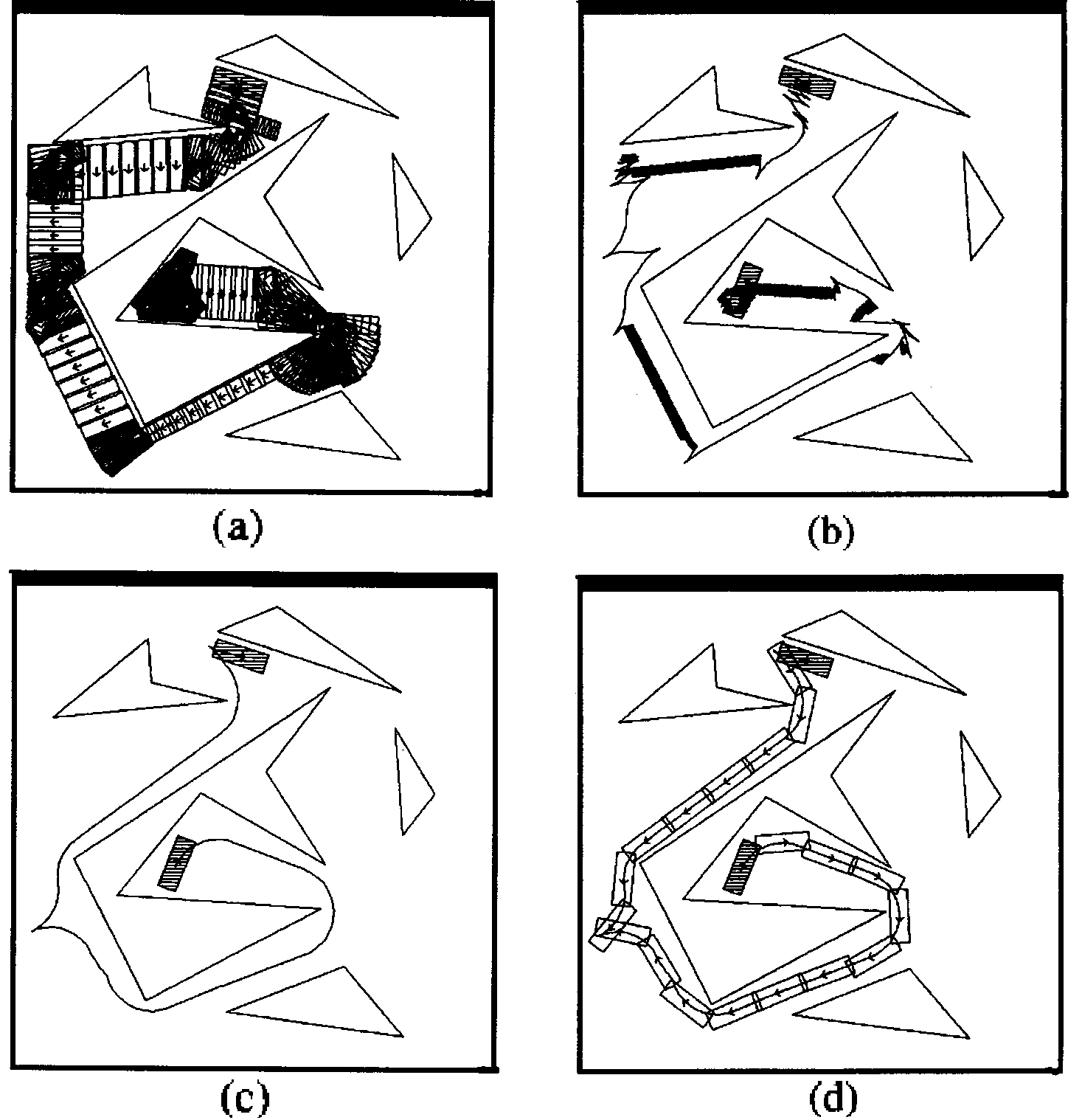} \qquad \includegraphics[width=0.45\linewidth]{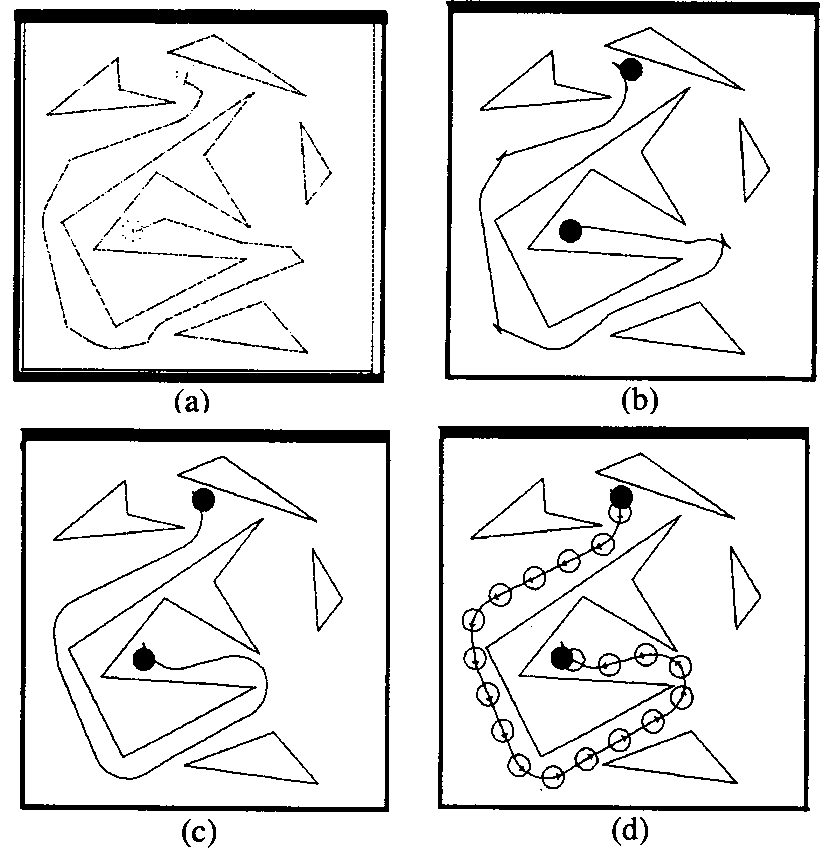}
\caption{Examples of trajectories generated by geometric planners from~\cite{laum94}, left pictures are made for a polygonal robot, right pictures for a circular robot: (a) holonomic path, (b) feasible path, (c) after optimization, (d) result}
\label{fig:laum94}
\end{figure}

In 1994, J.-P. Laumond, P. Jacobs, M. Ta\"{\i}x and R. Murray~\cite{laum94} developed a fast and accurate planner for a car-like model in presence of obstacles, based on a recursive subdivision of a collision free path generated by a lower-level geometric planner that ignores the motion constraints. The resultant trajectory is optimized to give a path that has near-minimal length in its homotopy class. The model of a front-wheel-drive car with a constraint on the turning radius can be reduced to system~(\ref{sys}) with the constraints 
\begin{align}
|u_2(t)| \leq |u_1(t)| \leq \frac{1}{R}. \label{constraint3}
\end{align}
The shortest paths for system~(\ref{sys}) satisfying inequalities~(\ref{constraint3}) are the same as for Reeds--Sheep's problem (see the proof in~\cite{sussmann}). 

The algorithm of the planner consists of three steps~\cite{laum94}: 
\begin{itemize}
\item[(a)] To plan a collision-free path for the geometric system, i.e., without taking into account differential system~(\ref{sys}). If such a path does not exist, neither does a feasible path. The experimental results present two implementations based on two different geometric planners. The first one runs for a polygonal robot (see left Fig.~\ref{fig:laum94} (a)), the second one for a disk based on the Voronoi diagram (see right Fig.~\ref{fig:laum94} (a)).
\item[(b)] To perform subdivision of the path until all endpoints are linked by a minimal-length path satisfying system~(\ref{sys}) and constraints~(\ref{constraint2}). The shortest path for Reeds-Shepp car is used in order to link intermediate configurations along the path, generated in the previous step (see~Fig.~\ref{fig:laum94}~(b)).
\item[(c)] To perform an "optimization" routine to remove extra maneuvers and reduce the length of the path (See Fig.~\ref{fig:laum94} (c)). The resulting paths are showed on Fig.~\ref{fig:laum94}~(d).
\end{itemize} 

Later in 1998, J.-P. Laumond presented notes~\cite{laum98} which describes a more general approach for the motion planning problem of a car and a car with one or two trailers. He mentioned several steering methods for step (b), including methods for nilpotent systems with piecewise and sinusoidal control. Back in those days, steering with optimal control was possible only for car-like systems, so, the only possibility for general systems was to call on numerical methods.

In 2013, H.~Chitsaz~\cite{chitsaz} considered the time-optimal control problem for system~(\ref{sysgen})--(\ref{vf}) with $l_r=0, l_t = 1$ and the constraints $(u_1, u_2) \in [-1,1]^2$. He revealed the structure of extremals through analytical integration of the system and adjoint equations. It was shown that Reeds-Shepp curves are time-optimal for system~(\ref{sysgen})--(\ref{vf}), and that the general case corresponds to a planar elastica (see Subsec.~\ref{secela}) connecting an arc of a circle to a straight segment. A complete characterization of optimal trajectories and control synthesis for the problem is still an open question.

This paper considers a sub-Riemannian problem for system~(\ref{sysgen})--(\ref{vf}), i.e., the problem with boundary conditions~(\ref{positions2}) and integral of sub-Riemannian length~(\ref{intsub}). Nowadays, computing optimal paths for such a problem is a challenging task. 

\section{Sub-Riemannian geometry}

Sub-Riemannian geometry considers optimal control problems with a system of linear in control differential equations and with the sub-Riemannian length minimization. For two-dimensional control, such a system has form~(\ref{sysgen}), where $X_1, X_2$ are smooth vector fields on manifold $M$ and $\dim(M) = n > 2$. The optimality criterion for two-dimensional control has form~(\ref{intsub}). System~(\ref{sysgen}) with~(\ref{positions}) and integral~(\ref{intsub}) describe an $n$-dimensional sub-Riemannian problem with 2-dimensional control.	

So far, only several sub-Riemannian problems have been studied completely, mostly for $n=3$. The simplest problems are nilpotent and are defined by the simplest structure of Lie algebras $\Lie (X_1,X_2)$. 

\subsection{Nilpotent approximation}
Control system~(\ref{sysgen}) is called nilpotent if Lie brackets of the corresponding vector fields $X_1, X_2$ vanish at a given length. Such systems provide a non-linear approximation for sub-Riemannian problems. Linearization of control systems is usually used as a local approximation, but for system~(\ref{sysgen}) which is linear in control the linearization is a too rough approximation: since the dimension of control is less than the dimension of configuration space, the linearization can not be controllable. Nilpotent approximation preserves this important property of a control system.

The term "nilpotent approximation" is defined by A.A.~Agrachev and A.V.~Sarychev in~\cite{agr-sar}, also by H.~Hermes in~\cite{hermes}. G.~Laferriere and H.J.~Sussmann propose a method in the context of nilpotent systems in~\cite{suss2}. Invariant construction of nilpotent approximation is proposed in~\cite{agr}. The exact solution to a control problem with a nilpotent system gives an approximate solution to arbitrary system~(\ref{sysgen}) in a small neighborhood of the final point.

The method of nipotent approximation can be applied to control problem~(\ref{sysgen})--(\ref{vf}) with arbitrary values of $l_t, l_r$. Different solutions to the nilpotent problem give different approximations. This article considers the solution in terms of optimal control~(\ref{intsub}) and uses a general algorithm, given by Bellaiche in~\cite{bellaiche} and further specified for system~(\ref{sysgen}) in~\cite{masht}, for constructing the nilpotent approximation. It uses the Maclaurin expansion of vector fields $X_1, X_2$ in privileged coordinates in the very same way as in~\cite{masht} for a mobile robot with two trailers. 

\subsection{Nilpotent sub-Riemannian problem on the Engel group}
The nipotent approximation for a mobile robot with a trailer is unique and has the following canonical form:
\begin{align}
\dot{\tilde{q}} &= u_1 \tilde{X}_1 (\tilde{q})  + u_2 \tilde{X}_2 (\tilde{q}), \label{syse}\\
\tilde{q} &= (\tilde{x},\tilde{y},z,v)^T \in \tilde{M}=\R^4, \quad (u_1, u_2) \in \R^2, \label{spacee}\\
\tilde{X}_1 &= \vect{1  \\ 0 \\  - \frac{\tilde{y}}{2} \\ 0}, \qquad \tilde{X}_2 = \vect{0 \\ 1 \\ \frac{\tilde{x}}{2} \\ \frac {\tilde{x}^2+\tilde{y}^2}{2}}. \label{vfe}
\end{align}
Problem~(\ref{syse})--(\ref{vfe}), (\ref{intsub}) with arbitrary boundary conditions is called the nilpotent sub-Riemannian problem on the Engel group. Since the problem is invariant under left shifts on the Engel group (see the defenition in~\cite{mont}), we can assume that the initial point is the identity, i.e.,
\begin{align}
\tilde{q}(0) = \tilde{q}_0 = (0,0,0,0), \qquad \tilde{q}(t_1) = \tilde{q}_1 = (\tilde{x}_1,\tilde{y}_1,z_1,v_1), \label{conden}
\end{align}
where $\tilde{q}_0, \tilde{q}_1 \in \tilde{M}$.

The corresponding sub-Riemannian problem has been studied quite recently in a series of works~\cite{engel, engel_conj, engel_cut}. The problem of finding optimal synthesis in the general case $\tilde{x}_1 z_1 \neq 0$ is reduced to a system of three algebraic equations in elliptic functions and elliptic integrals. It seems impossible to analytically solve such equations, therefore, a software for computing optimal trajectories for the nilpotent sub-Riemannian problem on the Engel group is being developed in Wolfram Mathematica~\cite{wolfram}. This software has already been used to devise several algorithms for computing approximate paths close to optimal in terms of~(\ref{intsub}) for a mobile robot with a trailer (see Subsec.~\ref{park}).

\section{Path planning via nilpotent approximation}
The algorithm proposed by A. Bellaiche in~\cite{bellaiche} is used to control a mobile robot with a trailer. Let us describe this algorithm applied to the problem~(\ref{sysgen})--(\ref{positions2}), (\ref{intsub}). Commutators of the vector fields $X_1, X_2$ taken from~(\ref{vf}) have the following form:
\begin{align*}
X_3=[X_1, X_2] = \vect{\sin \theta \\ -\cos \theta \\ 0 \\ -\frac{l_r + l_t \cos \varphi}{l_t^2}}, \quad X_4 = [X_1, X_3] = \vect{0 \\ 0 \\ 0 \\ -\frac{l_t + l_r \cos \varphi} {l_t^3}}.
\end{align*}

New coordinates $\tilde{q}$ are defined in the following way:
\begin{align}
\tilde{q} = \Gamma^{-1}(q-q_1). \label{gamma}
\end{align}
where $\Gamma$ is the matrix composed from vectors $X_1, X_2, X_3, X_4$ and vector $q_1$ corresponds to the final point of the curve, see~(\ref{positions2}). With such changes, terminal point $q = q_1$ goes to the origin $\tilde{q} = \tilde{q}_0$ and initial point $q = q_0$ goes to point $\tilde{q} = \tilde{q}_1$.

The next step is to find optimal control $(u_1,u_2) = (\hat{u}_1, \hat{u}_2)$ for problem~(\ref{sysgen}), (\ref{spacee}), (\ref{vfe}) with the boundary conditions $\tilde{q}(0) = \tilde{q}_1, \tilde{q}(t_1) = \tilde{q}_0$ and functional~(\ref{intsub}). It is possible to transform the boundary conditions to~(\ref{conden}) by inverting resulting controls and time on controls.  

The controls $(\hat{u}_1, \hat{u}_2)$ should be applied to initial problem~(\ref{sysgen})--(\ref{positions2}), (\ref{intsub}). The corresponding curve 
\begin{align}
\hat{q}(t) = \int_0^t \Big(X_1\big(q(t)\big) \hat{u}_1(t) + X_2\big(q(t)\big) \hat{u}_2(t) \Big) d t
\end{align} 
can be found by numerical integration.  It gives an approximate solution to the problem close to optimal in terms of~(\ref{intsub}). 

\subsection{Reparking of trailer}\label{reparking}
\begin{figure}[ht]
\centering
\includegraphics[width=0.99\linewidth]{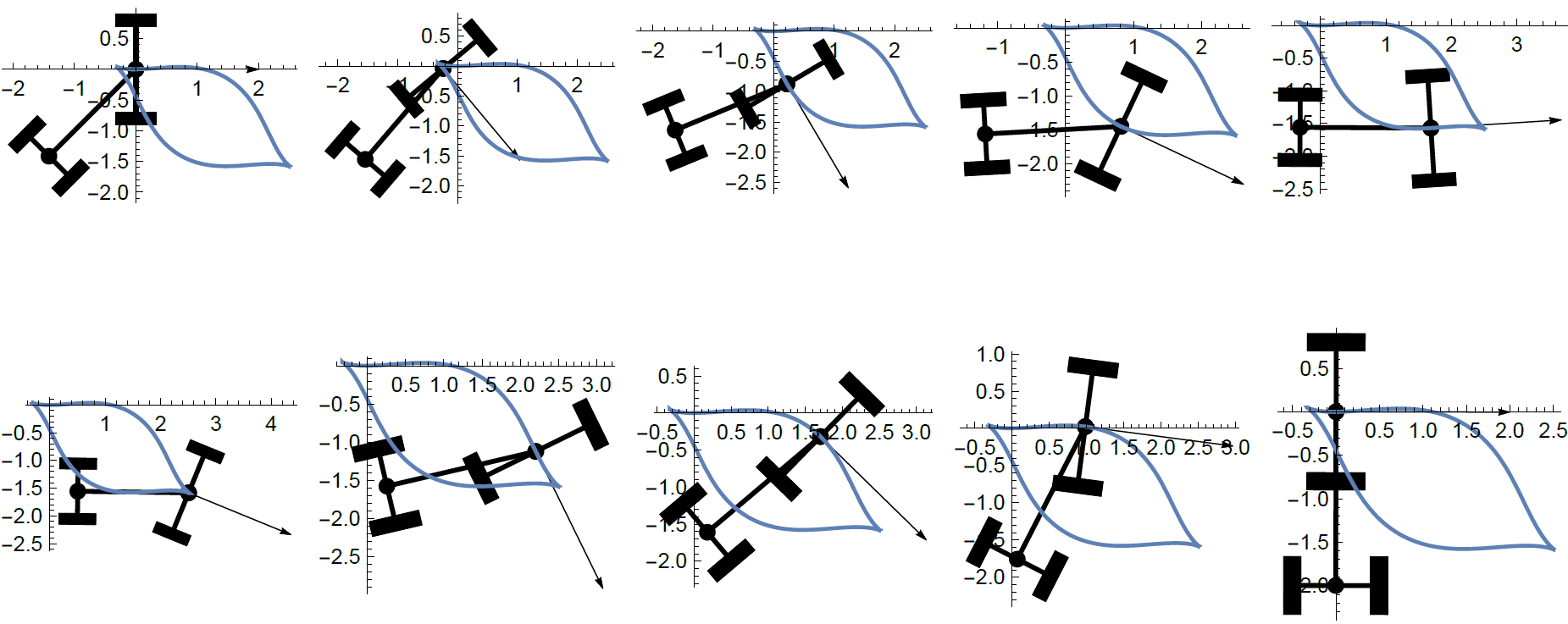}
\caption{Nilpotent approximation for reparking problem with parameters $(l_r,l_t) = (0, 2), (\varphi_0,\varphi_1) =(\frac{\pi}{4}, \frac{\pi}{2})$}
\label{fig:reparking1}
\end{figure}
Let us consider the particular case of problem (\ref{sysgen})--(\ref{positions2}),(\ref{intsub}):
\begin{align}
x_0 = y_0 = \theta_0 = x_1 = y_1 = \theta_1 = 0,
\end{align}
i.e. initial and final position of the car coincide one with each other, so the problem is to change the direction of the trailer. In this case, according to~(\ref{gamma}), the end point for the nilpotent problem on the Engel group is the following:
\begin{align}
\tilde{q}_1 = \left(0, 0, 0, \frac{l_t^3 (\varphi_1-\varphi_0)}{l_t + l_r \cos \varphi}\right). \label{v}
\end{align}

\begin{figure}[ht]
\centering
\includegraphics[width=0.99\linewidth]{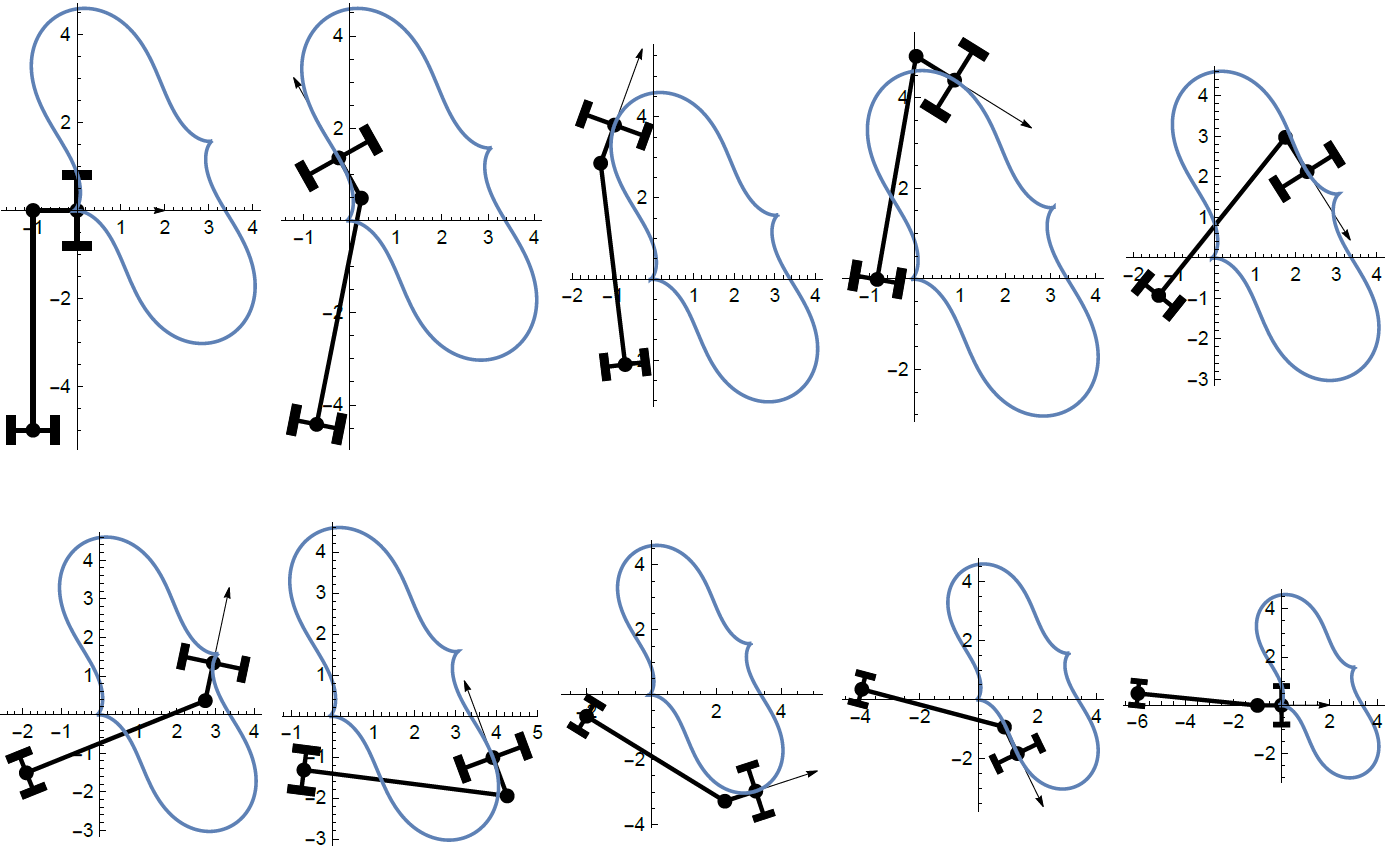}
\caption{Nilpotent approximation for reparking problem with parameters $(l_r,l_t) = (1, 5),  (\varphi_0, \varphi_1) = (\frac{\pi}{2}, -\frac{\pi}{3})$, $\alpha = 1$}
\label{fig:reparking2}
\end{figure}

For this case, there is an infinite number of optimal solutions which come to the same point $\tilde{q}_1$ in a figure-eight shape on the plane $(\tilde{x}, \tilde{y})$~\cite{engel_cut}. All of them can be used as the nilpotent approximation. The corresponding optimal controls are described by Jacobi elliptic functions:
\begin{align}
\hat{u}_1 (t) &= \sign (\sigma) 2 k_0 \sn \big( |\sigma| (\phi + \tcut - t) \big) \dn \big(|\sigma| (\phi + \tcut - t)\big), \label{cont1}\\
\hat{u}_2 (t) &= -\sign (\sigma) \Big(2 \dn^2 \big(|\sigma| (\phi + \tcut - t)\big) - 1\Big), \label{cont2}
\end{align}
where, in this case, $\tcut = \frac{4 K (k_0)}{|\sigma|}$ (see~\cite{engel_cut}) and $\phi \in [0,\tcut)$. On the other hand, using the parametrization of extremal trajectories~\cite{engel}, we get
\begin{align}
\tilde{q}_1 = \left(0,0,0,\frac{8 E(k_0)}{3 \sigma^3}\right) \qquad \Rightarrow  \qquad \sigma = \sqrt[3]{\frac{8 E(k_0)}{3 v_1}}. \label{sigma}
\end{align}

Numerical experiments show that a trajectory of system~(\ref{sysgen})--(\ref{vf}) starting from point $q(0) = (0,0,0, \varphi_0)$ with control~(\ref{cont1}),(\ref{cont2}) and arbitrary $\sigma, \phi$ always ends at point $q(\tcut) = (0,0,0, \tilde{\varphi}_1 (\sigma, \phi))$, where  $\tilde{\varphi}_1 (\sigma, \phi)$ is a function which can be calculated numerically. For fixed $\varphi_0,  \varphi_1$, the parameter $\sigma$ can be obtained from~(\ref{v}),(\ref{sigma}). A numerical algorithm is developed in order to find the value $\phi = \hat{\phi}$ which satisfies $\varepsilon = |\tilde{\varphi}_1(\sigma, \hat{\phi}) - \varphi_1| \to \min$. See an example of reparking obtained by this algorithm with accuracy $\varepsilon < 0.001$ in Fig.~\ref{fig:reparking1}. 

\begin{figure}[ht]
\centering
\includegraphics[width=0.99\linewidth]{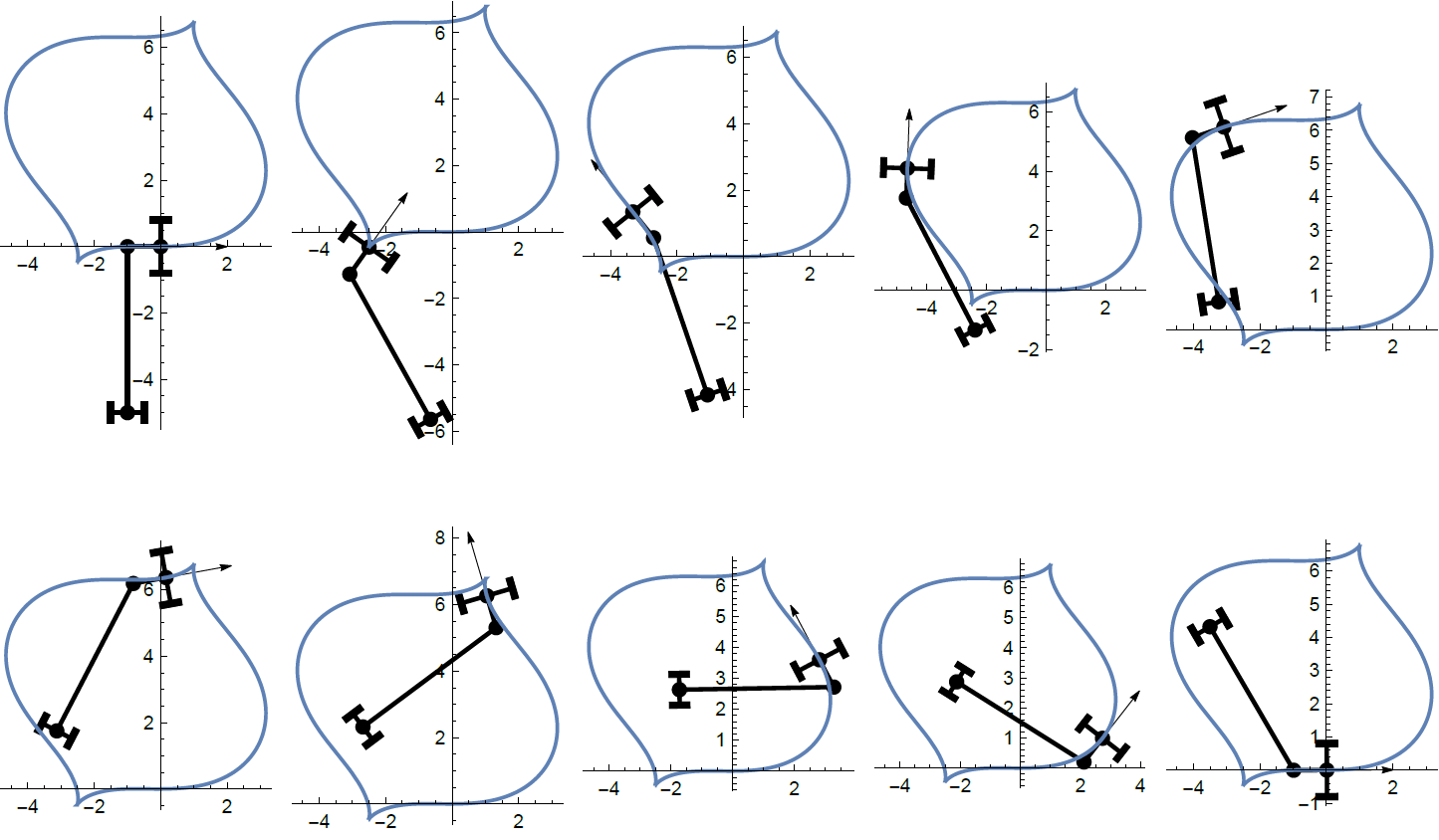}
\caption{Nilpotent approximation for reparking problem with parameters $(l_r,l_t) = (1, 5),  (\varphi_0, \varphi_1) = (\frac{\pi}{2}, -\frac{\pi}{3})$, $\alpha = 1.76113$}
\label{fig:reparking3}
\end{figure}

The algorithm was tested for different $\varphi_0, \varphi_1$ and $l_r, l_t$. Some cases provide rough accuracy $\varepsilon > 7/10$, see an example in Fig.~\ref{fig:reparking2}. Usually, such accuracy arises for distant points. One of the way to treat them is to use symmetry~(\ref{sym}). For the reparking problem it is possible to consider controls for models with $\delta_{\mu}(l_t)$, $\delta_{\mu}(l_r)$ and then translate them to controls for the model with $l_t, l_r$. Such a transformation also changes the functional integral to~(\ref{intsub2}) with $\alpha = \mu$. Therefore, using optimal solutions for the nilpotent sub-Riemannian problem, it is possible to obtain an approximate solution for each weighted integral~(\ref{intsub2}) with fixed $\alpha$. This consideration is equivalent to choosing scale on the plane. A similar trick has been used for sub-Riemannian problem on $\xSE(2)$ regarding cusp avoidance in image analysis applications~\cite{pde}. The example illustrated in Fig.~\ref{fig:reparking2} is solved with such consideration, see~Fig.~\ref{fig:reparking3}. 

\subsection{Algorithm for parking a mobile robot with a trailer} \label{park}
Let us consider problem~(\ref{sysgen})--(\ref{positions2}),(\ref{intsub}) with arbitrary $q_0, q_1$. The corresponding end point for the nilpotent problem on the Engel group according to~(\ref{gamma}) is expressed in the following way:
\begin{align}
\tilde{q}_1 =& \bigg((x_0 - x_1) \cos \theta_0 + (y_0 - y_1) \sin \theta_0, \quad \theta_0 - \theta_1, \nonumber \\
&\ (x_0 - x_1) \sin \theta_0 - (y_0 - y_1) \cos \theta_0, \quad  -\frac{l_t}{ l_t + l_r \cos \phi_0} \times \nonumber \\
&\ \times \Big(l_t^2 (\theta_0 - \theta_1 + \varphi_0-\varphi_1) + l_r l_t (\theta_0-\theta_1 ) \cos \varphi_0 + \nonumber \\
&\ + l_r \big((x_0-x_1 ) \sin \theta_0 - (y_0 - y_1) \cos \theta_0 \big)  + \nonumber \\
&\ + l_t \big((x_0-x_1 ) \sin (\varphi_0 + \theta_0) - (y_0 - y_1) \cos (\varphi_0 + \theta_0)\big) \Big)\bigg).
\end{align}

\begin{figure}[ht]
\centering
\includegraphics[width=0.99\linewidth]{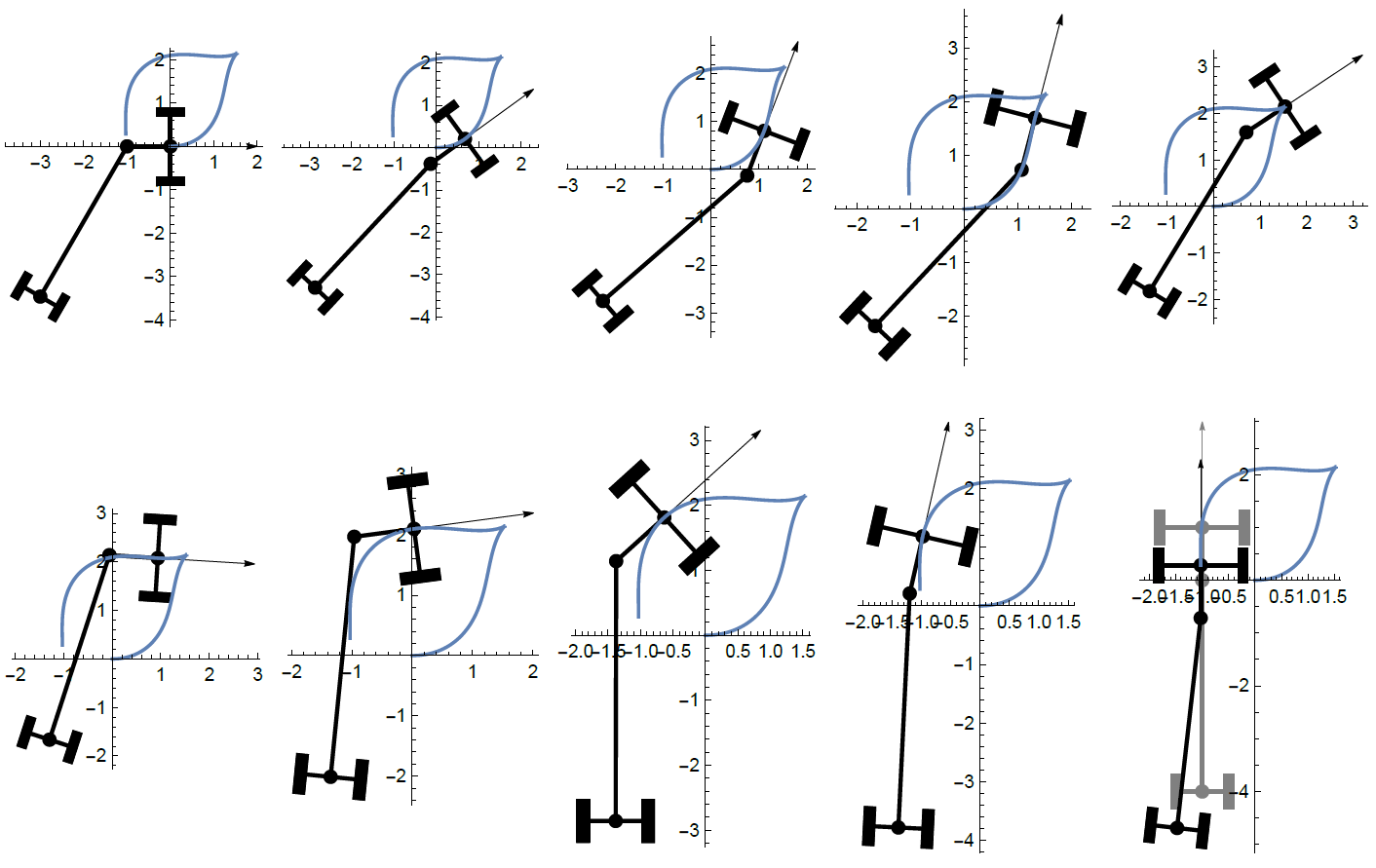}
\caption{Nilpotent approximation for parking problem with parameters $(l_r,l_t) = (1, 4),  q_0 = (0,0,0,\frac{\pi}{3}), q_1 = (-1, 1, \frac{\pi}{2}, 0)$}
\label{fig:parking1}
\end{figure}

If $\tilde{x}_1 z_1 \neq 0$, then there is a unique optimal trajectory for problem~(\ref{sysgen})--(\ref{positions2}),(\ref{intsub})~\cite{engel_cut}. Note that for case $\tilde{x}_1 z_1 = 0$ it is always possible to consider an arbitrarily close solution corresponding to a point $(\tilde{x}_1+\epsilon, \tilde{y}_1, z_1+\epsilon, v_1)$ with small enough $\epsilon$, which satisfies $(\tilde{x}_1+\epsilon) (z_1+\epsilon) \neq 0$. In order to find optimal controls for the end point with $\tilde{x}_1 z_1 \neq 0$, the corresponding system of algebraic equations in elliptic functions and elliptic integrals should be solved. Unfortunately, standard methods such as Newton's method, secant method, random search and grid search do not solve the system, therefore a hybrid method should be used instead. Such a method is developed by combining standard methods, and it finds a solution in most of the cases. A complete testing of the method is to be performed. We have reasonable grounds to believe that it is possible to develop an approach which can solve the system for any $\tilde{q}_1$. Thus, its detailed description is yet to be made in a future article.

\begin{figure}[ht]
\centering
\includegraphics[width=0.99\linewidth]{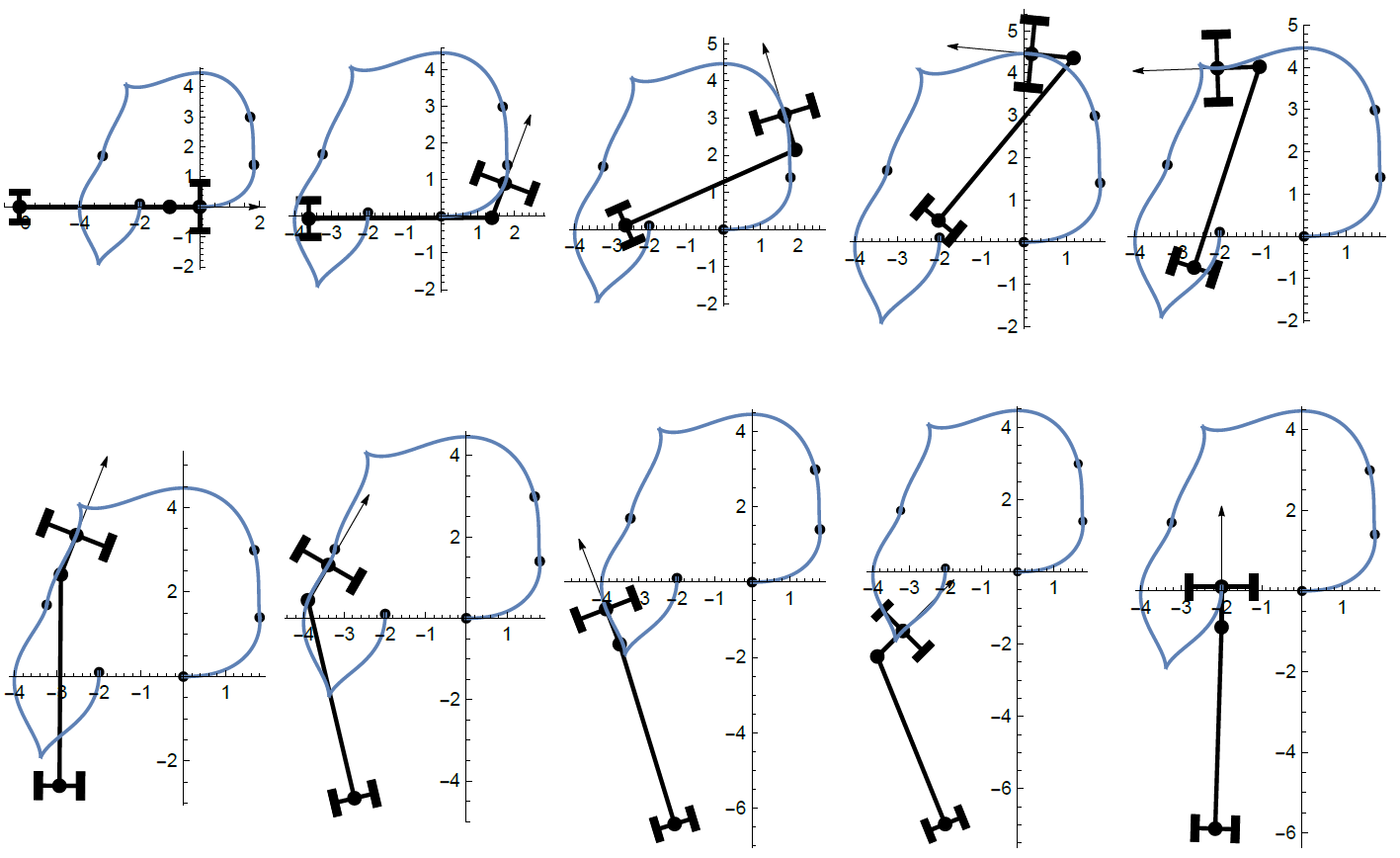}
\caption{Improved nilpotent approximation for parking problem with parameters $(l_r,l_t) = (1, 4),  q_0 = (0,0,0,\frac{\pi}{3}), q_1 = (-1, 1, \frac{\pi}{2}, 0)$}
\label{fig:parking2}
\end{figure}

For arbitrarily close initial and final configurations the nilpotent approximation gives solution $\hat{q}(t)$ close to optimal. However, such a solution usually does not connect the configurations. Starting from point $\hat{q}(0) = q_0$ it arrives at point $\hat{q}(\hat{t}_1) = \hat{q}_1 \neq q_1$. If the distance $\varepsilon = |\hat{q}_1 - q_1|$ is small enough, then the solution $\hat{q}(t), t\in [0, \hat{t}_1]$ can be used to solve the parking problem. Otherwise, there should be found a way to improve the obtained solution. Fig.~\ref{fig:parking1} shows an example with $\varepsilon > 0.72$. Here and below, gray images of a robot and a trailer correspond to desired positions.

One of the ways to improve the approximation is to consider the parking problem with the constraints $q(0)  = \hat{q}_0, q(t_1) = q_1$, where $\hat{q}_0 \neq q_0$ is an arbitrary point of curve $\hat{q}$. This improvement can be repeated several times until the resulting curve comes to a point close enough to $q_1$. Fig.~\ref{fig:parking2} shows an improved approximation of the example showed in Fig.~\ref{fig:parking1} with $\varepsilon < 0.022$. Black points correspond to points $\hat{q}_0$ where the nilpotent approximation has been recalculated.   

\begin{figure}[!ht]
\centering
\includegraphics[width=0.93\linewidth]{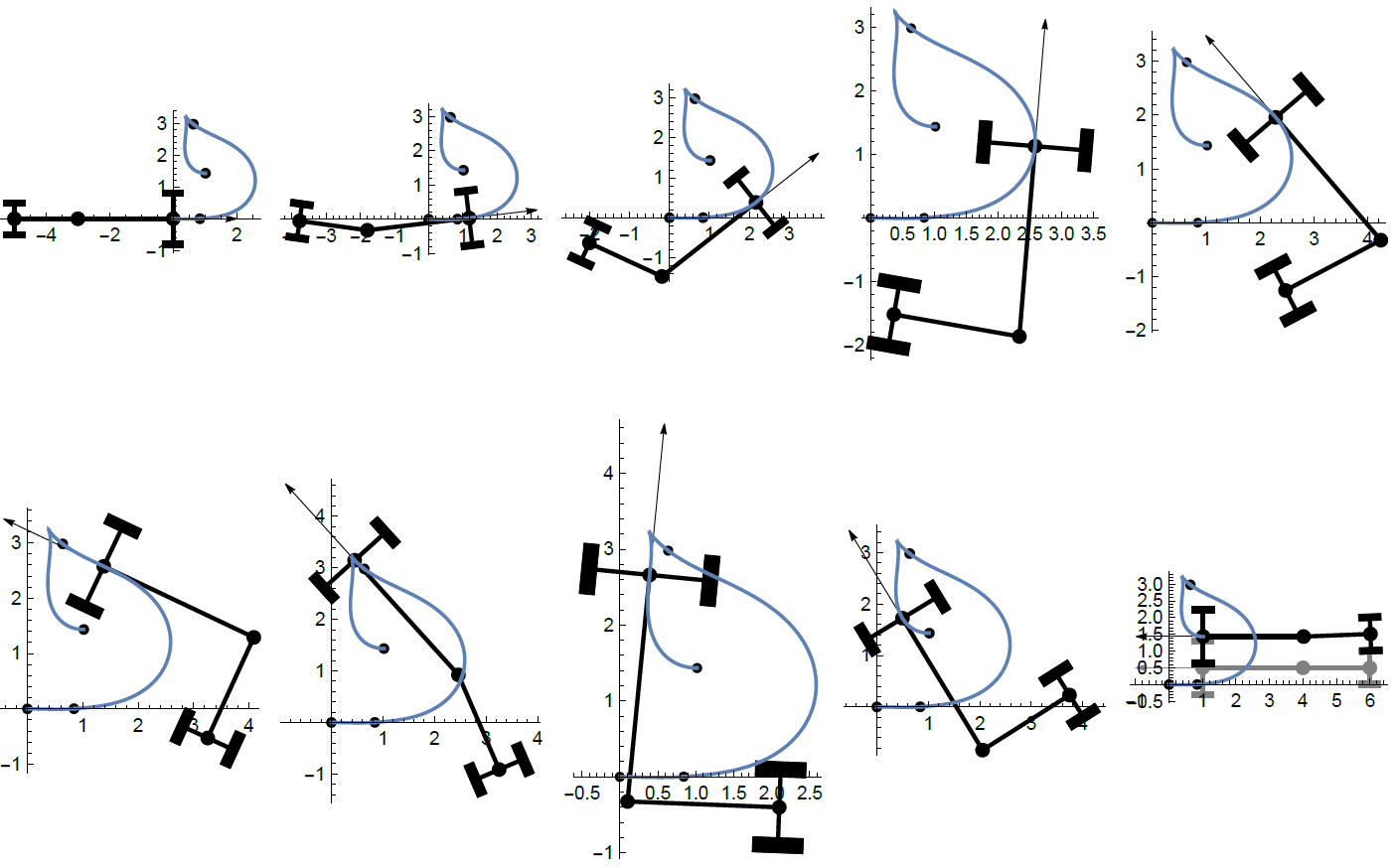}
\caption{Improved nilpotent approximation for a parking problem with parameters $(l_r,l_t) = (2, 3),  q_0 = (0,0,0,0), q_1 = (1, 1/2, \pi, 0)$}
\label{fig:parking3}
\end{figure}

Such an improvement is not always achieved with the desired precision, Fig.~\ref{fig:parking3}  illustrates an example with $\varepsilon \approx 0.935$ for improved approximation, where the initial approximation has precision $\varepsilon > 4.1$. In order to obtain a more accurate solution, a specified algorithm should be developed on the base of the proposed approach.

\begin{figure}[!ht]
\centering
\includegraphics[width=0.99\linewidth]{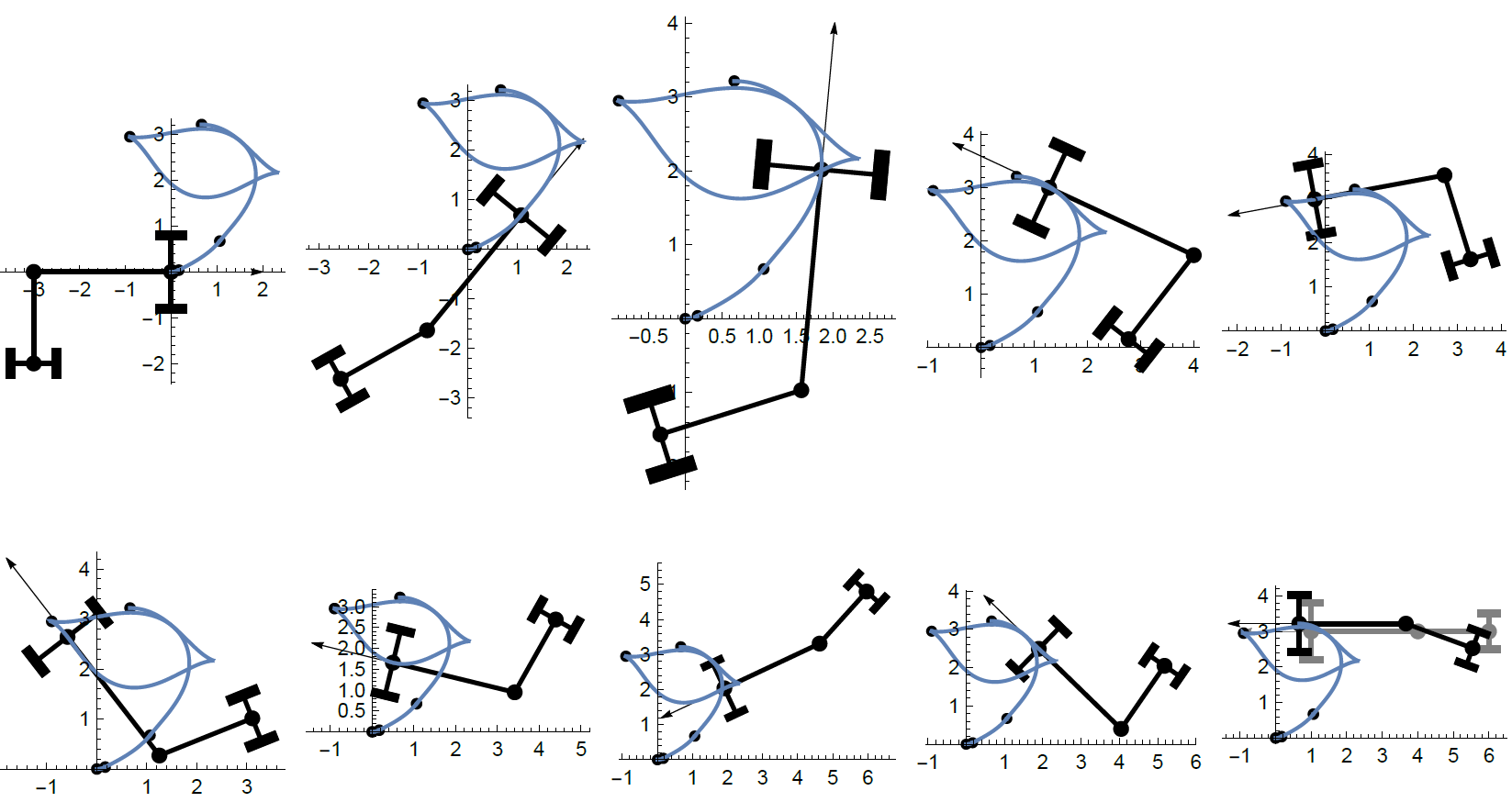}
\caption{Improved nilpotent approximation for a parking problem with parameters $(l_r,l_t) = (3, 2),  q_0 = (0,0,0,\frac{\pi}{2}), q_1 = (1, 3, \pi, 0)$}
\label{fig:parking4}
\end{figure}

\begin{figure}[!ht]
\centering
\includegraphics[width=0.99\linewidth]{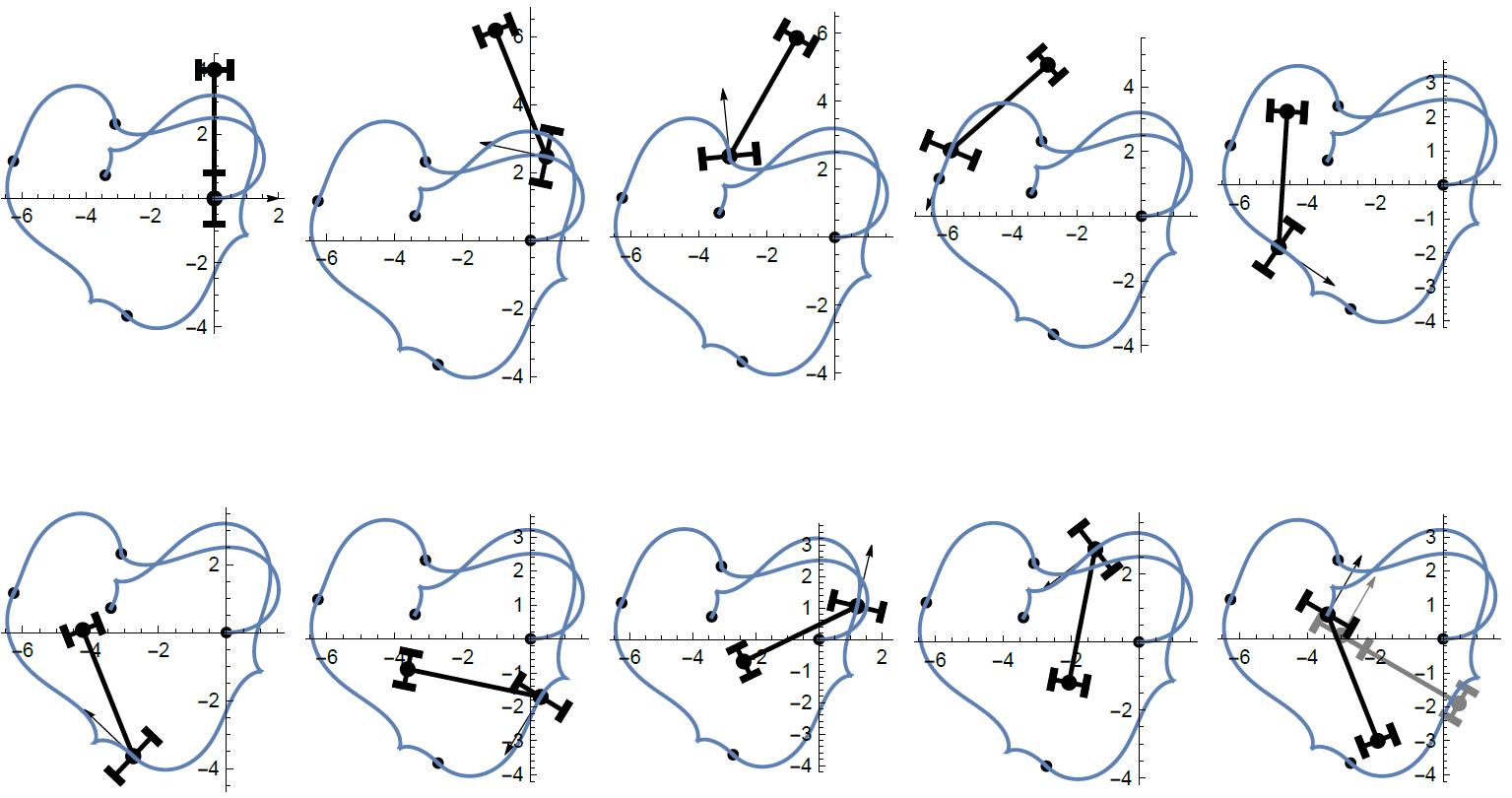}
\caption{Improved nilpotent approximation for a parking problem with parameters $(l_r,l_t) = (0, 4),  q_0 = (0,0,0,\frac{\pi}{2}), q_1 = (-3, \frac{1}{10}, \frac{\pi}{3}, \frac{\pi}{2})$}
\label{fig:parking5}
\end{figure}

More examples are shown in Fig.~\ref{fig:parking4} with $\varepsilon \approx 0.53$ (initial approximation has precision $\varepsilon > 6.11$) and in Fig.~\ref{fig:parking5} with $\varepsilon \approx 1$ (initial approximation has precision $\varepsilon > 3.8$). Note that, as the last example shows, there is no restrictions on the parameter $\varphi$.

\section*{Conclusion}
The article explores motion planning problems for a mobile robot with a trailer. This is quite a difficult task even without obstacles. The simplest case is given by the kinematic model of a car-like mobile robot, i.e. when the position of the trailer is not taken into account. The article also includes a brief overview of the existing methods of solving a motion planning problem for a mobile robot and a mobile robot with a trailer. One of them is based on the concept of nilpotent approximation. Different solutions, obtained from different classes of controls, are used to control a nilpotent system.  Optimal control, in sense of minimization of controls, corresponds to nilpotent sub-Riemannian problems.

Recently, the nilpotent sub-Riemannian problem on the Engel group was solved. This problem is given by a 4-dimensional control system with a 2-dimensional control and provides an approximate optimal solution to controlling a differential system for a mobile robot with a trailer. The algorithm of nilpotent approximation is applied for solving the problem of reparking a trailer. This algorithm is improved using a symmetry of the corresponding differential system for a mobile robot with a trailer. At the end of the article a new algorithm for parking a mobile robot with a trailer is presented. It will be applied for controlling a real mobile robot with a trailer.

The author acknowledges support by Russian Foundation for Basic Research, Project No. 16-31-00396-mol\_a.

\end{document}